\documentclass[12pt,reqno]{amsart}

\usepackage[latin1]{inputenc}
\usepackage{amsmath,amsfonts,amssymb,enumerate}
\usepackage{amsthm}
\usepackage{latexsym, mathtools,verbatim,graphicx}
\usepackage{amscd,hyperref}

\theoremstyle{plain}
\newtheorem{theorem}{Theorem}[section]
\newtheorem{thm}[theorem]{Theorem}
\newtheorem{prob}[theorem]{Problem}

\theoremstyle{definition}

\newtheorem{defn}[theorem]{Definition}

\theoremstyle{remark}
\newtheorem{rem}[theorem]{Remark}

\newcommand{\R}{\mathbb{R}}
\newcommand{\C}{\mathbb{C}}

\newcommand{\T}{\mathbb{T}}
\newcommand{\D}{\mathbb{D}}
\newcommand{\N}{\mathbb{N}}


\begin{document}
\title[Entire inner]{No entire inner functions}

\author[Cobos]{Alberto Cobos}
\address{Departement Wiskunde, KU Leuven, Celestijnenlaan 200B, B-3001 Leuven (Heverlee), Belgium.}\email{albertocobosrabano@gmail.com}
\author[Seco]{Daniel Seco}
\address{Instituto de Ciencias Matem\'aticas, Calle Nicol\'as Cabrera, UAM, 28049 Madrid,
Spain.} \email{dsf$\underline{\,\,\,}$cm@yahoo.es}

\date{\today}

\begin{abstract}
We study generalized inner functions on a large family of Reproducing Kernel Hilbert Spaces. We show that the only inner functions that are entire are the normalized monomials.
\end{abstract}

\maketitle

\section{Introduction}\label{Intro}

The Hardy space, $H^2$, is the space of functions $f$ holomorphic on the unit disc $\D$, with Taylor series $f(z) = \sum_{k=0}^\infty  a_k z^k$ satisfying \[\|f\|^2\coloneqq \sum_{k=0}^\infty |a_k|^2 < \infty.\]
A relevant property of $H^2$ functions is their well-defined boundary values at almost every point of the unit circle $\T$, defining a function $f^*$ whose $L^2(\T)$ norm equals the norm of $f$. A predominant role in $H^2$ is played by inner functions: $f \in H^2$  is \emph{inner} if $|f(z)| \leq 1$ for $z \in \D$ and \[|f^*(e^{i\theta})| = 1 \quad a.e.\ \theta \in [0, 2\pi).\]
This is equivalent to the function $f$ satisfying for all $k \in \N$
\begin{equation}\label{innerH2}
\left< z^k f, f \right> = \delta_k (0),
\end{equation}
where $\delta_k(0)$ is $1$ if $k=0$ and $0$, otherwise.  For more on $H^2$, see \cite{Gar}.

The following fact  is our starting observation:
\begin{rem}\label{rem1}
There are inner functions holomorphic on arbitrarily large discs, but if they are entire, they must be $\lambda z^k$, with $\lambda \in \T$. 
\end{rem}

Condition \eqref{innerH2} allows for generalization to other spaces:
Let $H$ be a Hilbert space of holomorphic functions over $\D$ with inner product $\left< \cdot , \cdot \right>_H$, and $f \in H$. We say that $f$ is \emph{$H$-inner} if for all $k \in \N$ it satisfies
\begin{equation}\label{innerH3}
\left< z^k f, f \right>_H = \delta_k (0).
\end{equation}
This goes back to Richter's work on invariant subspaces in the 80s and an account of properties of such functions is given in \cite{Remarks}.

A priori \eqref{innerH3} could depend on the norm. Indeed, some results in \cite{Remarks} have proved norm-dependent. This is a typical feature of invariant subspaces, to which $H$-inner functions are closely related.

In this note, we deal with the problem of whether the natural analogue of Remark \ref{rem1} is true in more general spaces.  A special role on the theory of $H$-inner functions is played by the analogue of (finite) Blaschke products, (finite) \emph{Shapiro-Shields functions}, which we describe in Section 2. As we will see then, Shapiro-Shields functions are examples of $H^2_\omega$-inner functions and they may be taken holomorphic over arbitrarily large discs. On the other hand, it is simple to see that a polynomial which is $H^2_\omega$-inner must be a normalized monomial. This justified posing the following problem in \cite{Secoprobs}:
\begin{prob}\label{theprob}
Is there any weight $\omega$ for which some (non-monomial) entire function becomes $H^2_\omega$-inner?
\end{prob}

Our setting will be so-called Dirichlet-type spaces, where some pathological case arises: in some spaces it is not possible to choose a Shapiro-Shields function with a prescribed zero set (and multiplicities). Zeros that do not match what prescribed  are called \emph{extraneous}. Many Dirichlet-type spaces do not have this pathology, although some choices of equivalent norms in $D_\alpha$ have been shown to have extraneous zeros \cite{Nowak}. We do not know whether these pathologies are due to the choice of space or of equivalent norm, or whether they can be avoided. 

In Section 3, we do our contribution to Problem \ref{theprob}: \emph{in spaces without extraneous zeros, all entire inner functions are monomials}.

\section{Weighted Hardy spaces and Shapiro-Shields functions}

A particular class of inner functions are \emph{finite Blaschke products}: they are the only inner functions to be holomorphic \emph{across} the boundary $\T$ (this, we denote $Hol(\overline{\D})$). They consist of finite products of \emph{Blaschke factors}, this is, monomials and functions 
\[f(z)= \frac{z-a}{1 - \overline{a}z},\]
for some $a \in \D$.
On the other hand, the \emph{reproducing kernel} at the point $a \in \D$ is the only function with the property that for all $g \in H^2$, \[g(a) = \left< g, k_a \right>,\] and it is given by
\[k_a(z)= \frac{1}{1-\overline{a}z}.\]
A Blaschke factor may be expressed as a normalization of the function $k_a(z)-k_a(a)$. As described in \cite{Remarks}, the analogue of finite Blaschke products in other Hilbert spaces with reproducing kernels (RKHS) is based on this. We focus on the following family of RKHS:
\begin{defn}\label{spaces}
Let $\omega = (\omega_k)_{k \in \N} \subset (0, +\infty)$, with $\omega_0=1$, satisfying \begin{equation}\label{conditionweight}
\lim_{k \rightarrow \infty} \frac{\omega_k}{\omega_{k+1}} = 1.\end{equation}  
The \emph{weighted Hardy space} with weight $\omega$ is the space $H^2_\omega$ of holomorphic functions $f$ over $\D$  given by $f(z) = \sum_{k=0}^{\infty} a_k z^k$  with 
\begin{equation}\label{norm}
\|f\|_{\omega}^2 \coloneqq  \sum_{k = 0}^{\infty}  |a_k|^2 w_k < \infty.
\end{equation}
We also assume a doubling condition: the existence of a constant $C>0$ such that \begin{equation}\label{doubling} \sup_{n \in \N} \sup_{n \leq k \leq 2n} \frac{\omega_k}{\omega_n} \leq C.\end{equation}\end{defn}

Assumption \eqref{conditionweight} guarantees both forward and backward shifts are bounded, and $Hol(\overline{\D}) \subset H^2_\omega \subset Hol(\D)$. A special class is formed by the weights $\omega_k = (k+1)^\alpha$, for some $\alpha \in \R$ (called \emph{Dirichlet-type spaces}). Specially well-known are the 3 cases when $\alpha = -1,0,1$, which are respectively the Bergman ($A^2$), Hardy ($H^2$) and Dirichlet ($D$) classical spaces. For more on $A^2$ and $D$, see \cite{DuS,EFKMR}. The role of the doubling condition is to ensure that spaces with bounded point evaluation functionals become multiplicative algebras (see Remark 3.2(II) below).

As hinted, $D_\alpha$ spaces are RKHS and the reproducing kernel at $a \in \D$ is given by \begin{equation}\label{eqn101}
k_a(z)=  \sum_{k=0}^\infty \frac{\overline{a}^kz^k}{\omega_k}. \end{equation}
In the notation below, whenever a function's zeros are denoted with some of them appearing more than once, we interpret that the number of appearances represents the multiplicity of the zero of the function at that point. We are now ready to define one of the main concepts we use:
\begin{defn}\label{ShaShi}
Let $H^2_\omega$ be a fixed weighted Hardy space, $k_a$ be the reproducing kernel of $H^2_\omega$ at any point $a \in \D$, and $Z=(z_k)_{k=1}^t \subset \D \backslash \{0\}$, where each $z_k \in Z$ has multiplicity one. The \emph{Shapiro-Shields function} corresponding to $Z$ (in $H^2_\omega$) is the function $h_Z$ defined as \[h_Z(z) \coloneqq \frac{f_Z(z)}{\|f_Z\|},\] where $f_Z(z)$ is given by the $(t+1) \times (t+1)$ determinant
\begin{equation}\label{eqn201}
f_Z(z) =
\begin{vmatrix}
1 	& 1\ldots 1 \\
(k_{z_j}(z))_{j=1}^t 	& (k_{z_j}(z_k))_{j,k=1}^t
\end{vmatrix}.\end{equation}
If the points are repeated or equal to 0, the corresponding Shapiro-Shields function is the limit of Shapiro-Shields functions for any sequence of $t$-tuples of non-repeated points, converging to the prescribed zero set in $\C^t$. Whenever the reproducing kernel exists on points of the boundary $\T$, one can also consider equation \eqref{eqn201} to define the Shapiro-Shields function for a zero set that contains some points of the boundary.
\end{defn}

Shapiro-Shields functions originate from \cite{ShaShi} and they have been studied in detail in the article \cite{Remarks}, among others. Some of their basic properties are the following:
\begin{rem}\label{remshashi} Let $Z \subset \D$ and $h_Z$ be the Shapiro-Shields function corresponding to $Z$ in the $H^2_\omega$ space.
\begin{itemize}
\item[(1)] $h_Z$ is $H^2_\omega$-inner.
\item[(2)] $h_Z(z_j)= 0$ for all $z_j \in Z$.
\item[(3)] $h_Z$ is holomorphic across the boundary but singular at $1/\overline{z_j}$ for any $z_j \in Z \backslash \{0\}$.
\end{itemize}
\end{rem}
Property (2) also holds for points of $Z$ on $\T$ and if we prescribe multiplicity $k$ then also the first $k-1$ derivatives will be 0 at the prescribed point. However, in the case of a $D_\alpha$ space, for $\alpha \in (1,2]$ it may very well happen that multiplicity $k$ at $z_j \in \T$ does not mean that $|h_Z(z)| \leq C|(z-z_j)^k|$ for any $C \in \R$. Indeed, it is simple to see this for $k=1$, $z_1=1$, since the derivative of $k_1$ at $1$ is unbounded.
This is one of the pathologies that we will need to avoid for our arguments to work. In general, we will want to avoid multiplicities of the zeros of Shapiro-Shields functions to be different from the prescribed ones: 
\begin{defn}\label{extraneous}
Let $h_Z$ be the Shapiro-Shields function given by the prescribed zero set $Z \subset \overline{\D}$ in $H^2_\omega$, let $z_0 \in \overline{\D}$ and let $k$ be the number of times $z_0$ appears on $Z$ (possibly, 0). We say that  $h_Z$ is \emph{regular}  at the point $z_0 \in \overline{\D}$ if, as $z$ tends to $z_0$, $|h_Z(z)|$ is comparable to $|z-z_0|^k$. Otherwise, we say that $z_0$ is an \emph{extraneous zero} of $h_Z$.
If $H^2_\omega$ contains no Shapiro-Shields functions with extraneous zeros, we will say that $H^2_\omega$ \emph{has no extraneous zeros}.
\end{defn}

In view of the paper \cite{Nowak}, it seems like the extraneous zero phenomenon may happen in rather simple spaces even with points inside the disc, whereas in principle, the example of $D_\alpha$ for $\alpha \in (1,2]$ happens only for points on the boundary.
Describing the spaces in which extraneous zeros appear seems like too ambitious of a goal for our purposes, and therefore we concentrate on the case of regular zeros.

\section{Invariant subspaces and Shapiro-Shields functions}

We are finally ready to state our main result:

\begin{thm}\label{thethm}
Let $H^2_\omega$ have no extraneous zeros. Then the only entire functions which are $H^2_\omega$-inner are the normalized monomials.
\end{thm}

Before the proof, let us recall a few facts about $H^2_\omega$-inner functions and point evaluations in $H^2_\omega$ (see \cite{Remarks} and \cite{FMS}, respectively). For $f$ in a space $H^2_\omega$, we denote $[f]_\omega$ the smallest closed subspace of $H^2_\omega$ which is invariant under multiplication by $z$ and contains $f$. Also, $H^\infty$ is the space of bounded analytic functions on $\D$.
\begin{rem}\label{rem301}
Let $\omega$ be as in Definition \ref{spaces}, fixed.
\begin{itemize}
\item[(I)] Let $g$ be $H^2_\omega$-inner, then the orthogonal projection of the function 1 onto $[g]_\omega$ is $g\overline{g(0)}$. In particular, if $g_1, g_2$ are both $H^2_\omega$ and $[g_1]_\omega=[g_2]_\omega$ then $g_1 = c g_2$ for some $c \in \T$.
\item[(II)] Because of the doubling condition \eqref{doubling}, $H^2_\omega$ is a multiplicative algebra if and only if evaluation at a given point of the boundary is a bounded functional which happens if and only if \[\sum_{k=0}^{\infty} \frac{1}{\omega_k} < \infty.\] Indeed, from Proposition 32 in \cite{Shields}, a space with finite $\sum 1/\omega_k$ is an algebra provided that 
\[\sup_{n\in \N} \sum_{k=0}^n \frac{\omega_n}{\omega_k \omega_{n-k}} < \infty,\]
but under the doubling condition the above quantity is controlled by \[2C \sum_{k=0}^{n/2} 1/\omega_k < \infty.\] Without the doubling condition, additional difficulties would have appeared (see \cite{AHMCR,Salas}).
\item[(III)] Suppose that $H^2_\omega$ is a multiplicative algebra, $g_1, g_2 \in H^2_\omega$ and $g_1/g_2 \in H^\infty$. Then $g_1/g_2 \in H^2_\omega$ and $g_1 \in [g_2]_\omega$.
\end{itemize}
\end{rem}

Now we have all the ingredients we need.
\begin{proof}
Let $\omega$ be fixed. Assume then $f$ to be an entire $H^2_\omega$-inner function, with a zero set in $\overline{\D}$ given by $Z$ (with the multiplicities prescribed by $Z$). Since $f$ is entire, it can only have a finite number of points in $Z$, each with finite multiplicity. 

Suppose firstly that evaluation at points of $\T$ are \emph{not} bounded functionals in $H^2_\omega$ and that $f(0) \neq 0$.
Define then $h_{Z_0}$ as the Shapiro-Shields function defined in $H^2_\omega$ by the zero set $Z_0= Z \cap \D$.  By our hypothesis, $h_{Z_0}/f$ and $f/h_{Z_0}$ are both holomorphic across the boundary and hence they are multipliers in $H^2_\omega$. This means that $h_{Z_0} \in [f]_\omega$ and that $f \in [h_{Z_0}]_\omega$, which, in turn, implies that \begin{equation}\label{eqn302} [f]_\omega = [h_{Z_0}]_\omega. \end{equation}
Since $f$ is $H^2_\omega$-inner, \eqref{eqn302} together with Remark \ref{rem301} (I) imply that $f= \lambda h_{Z_0}$ for some $\lambda \in \T$, but by Remark \ref{remshashi} (3), $h_{Z_0}$ is not entire, and so we arrive to a contradiction unless $Z_0 = \emptyset$. This means that $f$ must be a unimodular constant.
If $0 \in Z$, the same argument can be worked out with limits, but the conclusion will be that $f$ is $\lambda z^s$ for some $s \in \N$, $\lambda \in \T$.

Finally, suppose we are in the case that evaluation at some points of the boundary are bounded functionals. Bearing in mind Remark \ref{rem301} (II), all points of the boundary have this property, and the space $H^2_\omega$ is an algebra. Given $f$ with zeros in $\overline{\D}$ given by $Z$, we construct $h_Z$. Since $h_Z/f$ and $f/h_Z$ are bounded functions, they are also in the space, $h_Z \in [f]_\omega$ and $f \in [h_Z]_\omega$. We may conclude that $[f]_\omega=[h_Z]_\omega$. However, $h_Z$ is not entire unless a monomial, and so $f$ must be a monomial itself.
\end{proof}

\noindent\textbf{Acknowledgements.} The authors acknowledge financial support by
the Severo Ochoa Programme for Centers of Excellence in R\&D
(SEV-2015-0554) at ICMAT. The second author is also grateful for support from the
Spanish Ministry of Economy and Competitiveness, through grant
MTM2016-77710-P. We also thank I. Efraimidis and M. Hartz for careful reading and valuable comments.

\end{document}